\documentclass[11pt]{article}
\usepackage[matrix,arrow,curve,cmtip]{xy}
\usepackage{amssymb}
\usepackage{theorem}
\usepackage{latexsym}
\usepackage[a4paper,width=140mm,height=220mm]{geometry}
\usepackage{titlesec}

\titleformat{\section}[hang]%
{\bfseries\filcenter\large}{\thesection.}{1ex}{}%

\def\eqref#1{(\ref{#1})}

\def\to{\mbox{$\xymatrix@1@C=5mm{\ar@{->}[r]&}$}}
\def\tto{\mbox{$\xymatrix@1@C=5mm{\ar@{=>}[r]&}$}}
\def\halfcirc{\begin{picture}(0,0)\put(0,3){\oval(2,2)[l]}\end{picture}}
\def\incl{\mbox{$\xymatrix@1@C=5mm{\ar@{->}[r]|<{\halfcirc}&}$}}
\def\bkar{\mbox{$\xymatrix@1@C=5mm{\ar@{->}[l]&}$}}
\def\distsign{\begin{picture}(0,0)\put(0,0){\circle{4}}\end{picture}}
\def\dist{\mbox{$\xymatrix@1@C=5mm{\ar@{->}[r]|{\distsign}&}$}}
\def\bkdist{\mbox{$\xymatrix@1@C=5mm{\ar@{->}[l]|{\distsign}&}$}}
\def\biar{\mbox{$\xymatrix@1@C=5mm{\ar@<1.5mm>[r]\ar@<-0.5mm>[r]&}$}}
\def\bidist{\mbox{$\xymatrix@1@C=5mm{\ar@<1.5mm>[r]|{\distsign}\ar@<-0.5mm>[r]|{\distsign}&}$}}
\def\adjar{\mbox{$\xymatrix@1@C=5mm{\ar@<1.5mm>@{<-}[r]\ar@<-0.5mm>[r]&}$}}
\def\adjdist{\mbox{$\xymatrix@1@C=5mm{\ar@<1.5mm>@{<-}[r]|{\distsign}\ar@<-0.5mm>[r]|{\distsign}&}$}}
\def\iso{\mbox{$\xymatrix@1@C=6mm{\ar@{->}[r]^{\sim}&}$}}
\def\doubiso{\mbox{$\xymatrix@1@C=6mm{\ar@{<->}[r]^{\sim}&}$}}
\def\doubar{\mbox{$\xymatrix@1@C=6mm{\ar@{<->}[r]&}$}}
\def\endoar#1#2{\mbox{\xymatrix{{#1}\ar@(u,r)^{#2}}}}

\newtheorem{theorem}{Theorem}[section]
\newtheorem{lemma}[theorem]{Lemma}
\newtheorem{definition}[theorem]{Definition} 
\newtheorem{proposition}[theorem]{Proposition}
\newtheorem{corollary}[theorem]{Corollary}
{\theorembodyfont{\upshape}\newtheorem{example}[theorem]{Example}}
{\theorembodyfont{\upshape}}
\newcommand{\proof}{\noindent {\it Proof\ }: }
\newcommand{\sketchofproof}{\noindent {\it Sketch of proof\ }: }
\def\endofproof{$\mbox{ }\hfill\Box$\par\vspace{1.8mm}\noindent}

\def\Idm{{\sf Idm}}
\def\Mod{{\sf Mod}}
\def\Cocont{{\sf Cocont}}
\def\Ord{{\sf Ord}}

\def\Sh{{\sf Sh}}

\def\:{\colon}
\def\1{{\bf 1}}
\def\impl{\Longrightarrow}
\def\2{{\bf 2}}

\def\QUANT{{\sf QUANT}}
\def\cc{_{\sf cc}}
\def\op{^{\sf op}}

\def\Sup{{\sf Sup}}

\def\Dist{{\sf Dist}}
\def\Cat{{\sf Cat}}
\def\Map{{\sf Map}}
\def\id{{\sf id}}
\def\Id{{\sf Id}}
\def\Q{{\cal Q}}
\def\R{{\cal R}}
\def\A{{\cal A}}
\def\B{{\cal B}}
\def\C{{\cal C}}

\def\P{{\cal P}}

\def\Z{{\cal Z}}

\def\lim{\mathop{\rm lim}}
\def\bbA{\mathbb{A}}
\def\bbB{\mathbb{B}}
\def\bbC{\mathbb{C}}

\def\tensor{\otimes}

\title{$\Q$-modules are $\Q$-suplattices}
\author{Isar Stubbe\footnote{Postdoctoral Fellow of the Research Foundation Flanders (FWO), Department of Mathematics and Computer Science, University of Antwerp, Middelheimlaan 1, 2020 Antwerpen, Belgium, {\tt isar.stubbe@ua.ac.be}}}
\date{January 2007\footnote{Published in {\em Theory~Appl.~Categ.} on July 15, 2007. Please refer to the published version, which is freely available on {\tt http://www.tac.mta.ca/tac/}.}}

\begin{document}

\maketitle

\begin{abstract}
It is well known that the internal suplattices in the topos of sheaves on a locale are precisely the modules on that locale. Using enriched category theory and a lemma on KZ doctrines we prove (the generalization of) this fact in the case of ordered sheaves on a small quantaloid. Comparing module-equivalence with sheaf-equivalence for quantaloids and using the notion of centre of a quantaloid, we refine a result of F. Borceux and E. Vitale.
\end{abstract}

\section{Introduction}\label{A}
When studying topos theory one inevitably must study order theory too: if only because many advanced features of topos theory depend on order-theoretic arguments using the internal Heyting algebra structure of the subobject classifier in a topos, as C. J. Mikkelsen [1976] illustrates plainly. One of the results of [Mikkelsen, 1976] states that an ordered object in an elementary topos $\cal E$ is cocomplete, i.e.\ it is an internal suplattice, if and only if the ``principal downset embedding'' from that object to its powerobject has a left adjoint in $\Ord({\cal E})$. In the case of a localic topos, it turns out that the internal suplattices in $\Sh(\Omega)$ are precisely the $\Omega$-modules, and supmorphisms are just the module morphisms [Joyal and Tierney, 1984; Pitts, 1988].
\par
Now consider quantaloids (i.e.\ $\Sup$-enriched categories) as non-commutative, multi-typed generalization of locales. Using the theory of categories enriched in a quantaloid, and building further on results by B. Walters [1981] and F. Borceux and R. Cruciani [1998], I. Stubbe [2005b] proposed the notion of {\em ordered sheaf on a (small) quantaloid $\Q$} (or {\em $\Q$-order} for short): one of several equivalent ways of describing a $\Q$-order is to say that it is a Cauchy complete category enriched in the split-idempotent completion of $\Q$. There is thus a locally ordered category $\Ord(\Q)$ of $\Q$-orders and functors between them. If one puts $\Q$ to be the one-object suspension of a locale $\Omega$, then $\Ord(\Omega)$ is equivalent to $\Ord(\Sh(\Omega))$. (And if one puts $\Q$ to be the one-object suspension of the Lawvere reals $[0,\infty]$, then $\Ord([0,\infty])$ is equivalent to the category of Cauchy complete generalized metric spaces.)
\par
In this paper we shall explain how $\Mod(\Q)$, the quantaloid of $\Q$-modules, is the category of Eilenberg-Moore algebras for the KZ doctrine on $\Ord(\Q)$ that sends a $\Q$-order $\bbA$ to its free cocompletion $\P\bbA$. The proof of this fact is, altogether, quite straightforward: a lot of the hard work -- involving quantaloid-enriched categories -- has already been done elsewhere [Stubbe, 2005a, 2005b, 2006], so we basically only need a lemma on KZ doctrines to put the pieces of the puzzle together. Applied to a locale $\Omega$, and up to the equivalence of $\Ord(\Omega)$ with $\Ord(\Sh(\Omega))$, this KZ doctrine sends an ordered sheaf on $\Omega$ to (the sheaf of) its downclosed subobjects, so our more general theorem provides an independent proof of the fact that $\Omega$-modules are precisely the internal cocomplete objects of $\Ord(\Sh(\Omega))$. This then explains the title of this paper: even for a small quantaloid $\Q$, ``$\Q$-modules are $\Q$-suplattices''! We end the paper with a comment on the comparison of (small) quantaloids, their categories of ordered sheaves, their module categories, and their centres; thus we refine a result of F. Borceux and E. Vitale [1992].
\par
In some sense, this paper may be considered a {\em prequel} to [Stubbe, 2007]: we can now rightly say that the latter paper treats those $\Q$-suplattices (in their guise of cocomplete $\Q$-enriched categories) that are totally continuous (or supercontinuous, as some say). It is hoped that this will lead to a better understanding and further development of ``dynamic domains'', i.e.~``domains'' in $\Ord(\Q)$, so that applying general results to either $\Omega$ or $[0,\infty]$ then gives interesting results for ``constructive domains'' or ``metric domains''.
 
\section{Preliminaries}\label{B}
\subsection*{Quantales and quantaloids}
Let $\Sup$ denote the category of complete lattices and maps that preserve arbitrary suprema ({\em suplattices} and {\em supmorphisms}): it is symmetric monoidal closed for the usual tensor product. A {\em quantaloid} is a $\Sup$-enriched category; it is {\em small} when it has a set of objects; and a one-object quantaloid (most often thought of as a monoid in $\Sup$) is a {\em quantale}. A $\Sup$-functor between quantaloids is a {\em homomorphism}; $\QUANT$ denotes the (illegitimate) category of quantaloids and their homomorphisms. A standard reference on quantaloids is [Rosenthal, 1996].
\par
For a given quantaloid $\Q$ we write $\Idm(\Q)$ for the new quantaloid whose objects are the idempotent arrows in $\Q$, and in which an arrow from an idempotent $e\:A\to A$ to an idempotent $f\:B\to B$ is a $\Q$-arrow $b\:A\to B$ satisfying $b\circ e=b=f\circ b$. Composition in $\Idm(\Q)$ is done as in $\Q$, the identity in $\Idm(\Q)$ on some idempotent $e\:A\to A$ is $e$ itself, and the local order in $\Idm(\Q)$ is that of $\Q$. (Note that $\Idm(\Q)$ is small whenever $\Q$ is.) It is easy to verify that the quantaloid $\Idm(\Q)$ is the universal split-idempotent completion of $\Q$ in $\QUANT$, as the next lemma spells out.
\begin{lemma}\label{1}
If $\R$ is a quantaloid in which idempotents split, then, for any 
quantaloid $\Q$, the full embedding $i\:\Q\to\Idm(\Q)\:(f\:A\to B)\mapsto(f\:1_A\to 
1_B)$ determines an equivalence of quantaloids $-\circ i\:\QUANT(\Idm(\Q),\R)\to\QUANT(\Q,\R)$.
\end{lemma} 
\par
When $\Q$ is a small quantaloid, we write $\Mod(\Q)$ for $\QUANT(\Q\op,\Sup)$: the objects of this (large) quantaloid are called the {\em modules on $\Q$}. Since idempotents split in $\Sup$, it follows directly from \ref{1} that $\Mod(\Q)\simeq\Mod(\Idm(\Q))$.
\subsection*{Quantaloid-enriched categories}
A quantaloid is a bicategory and therefore it may serve itself as base for enrichment. The theory of quantaloid-enriched categories, functors and distributors is surveyed in [Stubbe, 2005a] where also the appropriate references are given. To make this paper reasonably self-contained we shall go through some basic notions here; we follow the notations of {\it op.~cit.}\ for easy cross reference.
\par
To avoid size issues we work with a small quantaloid $\Q$. A {\em $\Q$-category} $\bbA$ consists of a set $\bbA_0$ of `objects', a `type' function $t\:\bbA_0\to\Q_0$, and for any $a,a'\in\bbA_0$ a `hom-arrow' $\bbA(a',a)\:ta\to ta'$ in $\Q$; these data are required to satisfy 
$$\bbA(a'',a')\circ\bbA(a',a)\leq\bbA(a'',a)\mbox{\ \ \ and\ \ \ }1_{ta}\leq\bbA(a,a)$$ 
for all $a,a',a''\in\bbA_0$. A {\em functor} $F\:\bbA\to \bbB$ is a map $\bbA_0\to\bbB_0\:a\mapsto Fa$ that satisfies 
$$ta=t(Fa)\mbox{\ \ \ and\ \ \ }\bbA(a',a)\leq\bbB(Fa',Fa)$$ 
for all $a,a'\in\bbA_0$. For parallel functors $F,G\:\bbA\biar\bbB$ we put $F\leq G$ when $1_{ta}\leq\bbB(Fa,Ga)$ for every $a\in\bbA_0$. With the obvious composition and identities we obtain a locally ordered category $\Cat(\Q)$ of $\Q$-categories and functors.
\par
To give a {\em distributor} (or {\em module} or {\em profunctor}) $\Phi\:\bbA\dist\bbB$ between $\Q$-categories is to specify for any $a\in\bbA_0$, $b\in\bbB_0$, an arrow $\Phi(b,a)\:ta\to tb$ in $\Q$, such that 
$$\bbB(b,b')\circ\Phi(b',a)\leq\Phi(b,a)\mbox{\ \ \ and\ \ \ }\Phi(b,a')\circ\bbA(a',a)\leq\Phi(b,a)$$
for every $a,a'\in\bbA_0$, $b,b'\in\bbB_0$. Two distributors $\Phi\:\bbA\dist\bbB$, $\Psi\:\bbB\dist\bbC$ 
compose: we write $\Psi\tensor\Phi\:\bbA\dist\bbC$ for the 
distributor with elements
$$\Big(\Psi\tensor\Phi\Big)(c,a)=\bigvee_{b\in\bbB_0}\Psi(c,b)\circ\Phi(b,a).$$
The identity distributor on a $\Q$-category $\bbA$ is 
$\bbA\:\bbA\dist\bbA$ itself, i.e.\ the distributor with elements 
$\bbA(a',a)\:ta\to ta'$. We order parallel distributors $\Phi,\Phi'\:\bbA\bidist\bbB$ by ``elementwise comparison'': we define $\Phi\leq\Phi'$ to mean that $\Phi(b,a)\leq\Phi'(b,a)$ for every $a\in\bbA_0$, $b\in\bbB_0$. It is easily seen that $\Q$-categories and distributors form a quantaloid $\Dist(\Q)$.
\par
Every functor $F\:\bbA\to\bbB$ between $\Q$-categories represents an adjoint pair of distributors: 
\begin{itemize}
\itemsep=-4pt
\item the left adjoint $\bbB(-,F-)\:\bbA\dist\bbB$ 
has elements 
$\bbB(b,Fa)\:ta\to tb$,
\item the right adjoint $\bbB(F-,-)\:\bbB\dist\bbA$ 
has elements 
$\bbB(Fa,b)\:tb\to ta$.
\end{itemize}
The assignment $F\mapsto\bbB(-,F-)$ is a faithful 2-functor from 
$\Cat(\Q)$ to $\Dist(\Q)$; it gives rise to a rich theory of $\Q$-categories. We shall briefly explain two notions that play an essential r\^ole in the current work: cocompleteness and Cauchy completeness.
\subsection*{Cocompleteness and modules}
Given a distributor $\Phi\:\bbA\dist\bbB$ and a functor $F\:\bbB\to\bbC$, a functor $K\:\bbA\to\bbC$ is the {\em $\Phi$-weighted colimit of $F$} when it satisfies
$$\bbC(K-,-)=[\Phi,\bbC(F-,-)]$$
(and in that case it is essentially unique). The right hand side of this equation uses the adjunction between ordered sets
$$
\Dist(\Q)(\bbC,\bbA)\xymatrix@=15mm{\ar@<1mm>@/^2mm/[r]^{\Phi\tensor-}\ar@{}[r]|{\bot} & \ar@<1mm>@/^2mm/[l]^{[\Phi,-]}}\Dist(\Q)(\bbC,\bbB)$$
which surely exists since $\Dist(\Q)$ is a quantaloid. A functor is {\em cocontinuous} if it preserves all weighted colimits that happen to exist in its domain; and a $\Q$-category is {\em cocomplete} if it admits all weighted colimits. We write $\Cocont(\Q)$ for the subcategory of $\Cat(\Q)$ of cocomplete $\Q$-categories and cocontinuous functors. Much more can be found in [Stubbe, 2005a, sections 5 and 6].
\par
As stated in [Stubbe, 2006, 4.13] (but see also the references contained in that paper), $\Mod(\Q)$ and $\Cocont(\Q)$ are biequivalent locally ordered categories. Indeed, a $\Q$-module $M\:\Q\op\to\Sup$ determines a $\Q$-category $\bbA_M$: as object set take $(\bbA_M)_0=\uplus_{X\in\Q_0}MX$, then say that $tx=X$ precisely when $x\in MX$, and for $x\in MX$, $y\in MY$ let $\bbA_M(y,x)=\bigvee\{f\:X\to Y\mid Mf(y)\leq x\}$. A detailed analysis of why this $\bbA_M$ is cocomplete, and why every cocomplete $\Q$-category arises in this way, is precisely the subject of [Stubbe, 2006]; we shall not go into details here. 
\begin{corollary}\label{3}
For a small quantaloid $\Q$, $$\Cocont(\Q)\simeq\Mod(\Q)\simeq\Mod(\Idm(\Q))\simeq\Cocont(\Idm(\Q))$$
are biequivalent locally ordered categories.
\end{corollary}
\subsection*{Cauchy completeness and orders}
A $\Q$-category $\bbC$ is {\em Cauchy complete} if for any other $\Q$-category $\bbA$ the map
$$\Cat(\Q)(\bbA,\bbC)\to\Map(\Dist(\Q))(\bbA,\bbC)\:F\mapsto \bbC(-,F-)$$
is surjective, i.e.\ when any left adjoint distributor (also called {\em Cauchy distributor}) into $\bbC$ is represented by a functor. This is equivalent to the requirement that $\bbC$ admits any colimit weighted by a Cauchy distributor; and moreover such weighted colimits are {\em absolute} in the sense that they are preserved by any functor [Street, 1983]. We write $\Cat\cc(\Q)$ for the full subcategory of $\Cat(\Q)$ whose objects are the Cauchy complete $\Q$-categories. For more details we refer to [Stubbe, 2005a, section 7].
\par
Now we have everything ready to state an important definition from [Stubbe, 2005b].
\begin{definition}
For a small quantaloid $\Q$, we write $\Ord(\Q)$ for the locally ordered category $\Cat\cc(\Idm(\Q))$, and call its objects {\em ordered sheaves on $\Q$}, or simply {\em $\Q$-orders}.
\end{definition}
In fact, the definition of `$\Q$-order' in [Stubbe, 2005b, 5.1] is not quite this one: instead it is given in more ``elementary'' terms (avoiding the split-idempotent construction). But it is part of the investigations in that paper (more precisely in its section 6) that what we give here as definition is indeed equivalent to what was given there; and for the purposes of the current paper this ``structural'' definition is best. 
\par
The notion of $\Q$-order has the merit of generalizing two -- at first sight quite different -- mathematical structures: On the one hand, taking $\Q$ to be the (one-object suspension of) the Lawvere reals $[0,\infty]$, $\Ord([0,\infty])$ is the category of Cauchy complete generalized metric spaces [Lawvere, 1973]. On the other hand, taking $\Q$ to be the (one-object suspension of) a locale $\Omega$, $\Ord(\Omega)$ is the category of ordered objects in the topos $\Sh(\Omega)$ [Walters, 1981; Borceux and Cruciani, 1998]; obviously, this example inspired our terminology. For details we refer to [Stubbe, 2005b].

\section{Monadicity of $\Q$-modules over $\Q$-orders}\label{C}
Recall from [Kock, 1995] that a {\em Kock--Z\"oberlein (KZ) doctrine} on a locally ordered 2-category $\C$ is a 
monad $(T\:\C\to\C,\eta\:\Id_{\C}\tto T,\mu\:T\circ T\tto T)$ for 
which $T(\eta_C)\leq\eta_{TC}$ for any $C\in\C$. This precisely 
means that ``$T$-structures are adjoint to units''. Further on we shall encounter an instance of the following lemma.
\begin{lemma}\label{4}
Consider locally ordered 2-categories and 2-functors as in
$$
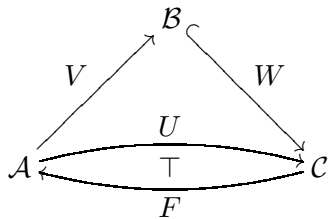
\xymatrix@=15mm{
   & \B\ar@{^{(}->}[dr]^W \\
\A\ar[ur]^V\ar@{}[rr]|{\top}\ar@/^3mm/[rr]^U & & 
\C\ar@/^3mm/[ll]^F}$$ with $W$ a local equivalence and $W\circ V=U$. Write $\eta\:\id_{\C}\tto U\circ F$ for the unit of the involved 
adjunction. Then 
\begin{enumerate}
\item $F\circ W\dashv V$ and its unit $\xi\:\Id_{\B}\tto V\circ(F\circ W)$ satisfies $\eta*\id_W = \id_V*\xi$, that is, $W(\xi_B)=\eta_{WB}$ for every $B\in\B$,
\newcounter{saveenum}
\setcounter{saveenum}{\value{enumi}}
\end{enumerate}
and writing $T=U\circ F\:\C\to\C$ and $S=V\circ (F\circ W)\:\B\to\B$, 
these monads satisfy
\begin{enumerate}
\setcounter{enumi}{\value{saveenum}}
\item $T\circ W = W\circ S$,
\item if $T$ is a KZ doctrine then 
\begin{enumerate}
\item also $S$ is a KZ doctrine,
\item $B\in\B$ is an $S$-algebra if and only if $WB$ is a $T$-algebra,
\item for $A\in\A$, $UA$ is a $T$-algebra if and only if $VA$ is an $S$-algebra,
\item if $\A\simeq\C^T$ then $\A\simeq\B^S$.
\end{enumerate}
\end{enumerate}
\end{lemma}
\proof To prove that $F\circ W\dashv V$, observe that for $B\in\B$ 
and $C\in\C$,
$$\xymatrix@=12mm{
\B(B,VC)\ar[d]^{\mbox{ apply }W} \\ \A(WB,WVC)\ar@{=}[d]^{\mbox{ use 
that }U=WV} \\ \A(WB,UC)\ar[d]^{\mbox{ use that }F\dashv U} \\ 
\C(FWB,C)}$$ are all equivalences (recall that $W$ is supposed to be 
a local equivalence). Putting $C=FWB$ in the above, and tracing the 
element $1_{FWB}$ through the equivalences, results in
$W(\xi_{B})=\eta_{WB}$. 
\par
The second part of the lemma is trivial.
\par
For the third part, suppose that $T(\eta_C)\leq\eta_{TC}$ for any 
$C\in\C$, then also
$$WS(\xi_B)=TW(\xi_B)=T(\eta_{WB})\leq\eta_{TWB}=\eta_{WSB}=W(\xi_{SB})$$
for every $B\in\B$; but $W$ is locally an equivalence, so 
$S(\xi_B)\leq\xi_{SB}$ as required to prove (a). Now, by the very 
nature of the algebras of KZ doctrines, $B\in\B$ is an $S$-algebra 
if and only if $\xi_B$ is a right adjoint in $\B$, which is the same 
as $W(\xi_B)=\eta_{WB}$ being a right adjoint in $\C$ because $W$ is 
locally an equivalence, and this in turn is just saying that $WB$ is 
a $T$-algebra. This proves (b), and (c) readily follows by putting 
$B=VA$ for an $A\in\A$, and using that $W\circ V=U$; so (d) becomes 
obvious.
\endofproof
In the rest of this section we let $\Q$ be a small quantaloid. It is a result from $\Q$-en\-riched category theory [Stubbe, 2005a, 6.11]
that $\Cocont(\Q)$, the locally ordered category of cocomplete $\Q$-categories and cocontinuous functors, is monadic over the locally ordered category $\Cat(\Q)$ of all categories and functors: the forgetful functor $\Cocont(\Q)\to\Cat(\Q)$ admits the presheaf construction as left adjoint, 
\begin{equation}\label{e1}
\Cocont(\Q)
\xymatrix{\ar@{}[r]|{\perp}\ar@/_2mm/@<-1mm>[r]_{\cal U} & 
\ar@/_2mm/@<-1mm>[l]_{\cal P}}\Cat(\Q).
\end{equation}
The unit of the adjunction is given by the Yoneda embeddings $Y_{\bbA}\:\bbA\to\P\bbA$; and a $\Q$-category $\bbA$ is in $\Cocont(\Q)$ if and only if $Y_{\bbA}\:\bbA\to\P\bbA$ admits a left adjoint in $\Cat(\Q)$, which 
is then the structure map of the algebra $\bbA$. In short, the monad induced by \eqref{e1} is a KZ-doctrine on $\Cat(\Q)$. 
\par
Cauchy complete $\Q$-categories can be characterized as those $\Q$-categories that admit all absolute colimits [Stubbe, 2005a, 7.2]. Knowing this it is clear that the forgetful $\Cocont(\Q)\to\Cat(\Q)$ factors over the full embedding $\Cat\cc(\Q)\to\Cat(\Q)$ of Cauchy complete $\Q$-categories into all $\Q$-categories. Applying \ref{4} to the adjunction in \eqref{e1} we thus obtain that the forgetful $\Cocont(\Q)\to\Cat\cc(\Q)$ has a left adjoint, namely (the restriction of) the presheaf construction, and moreover $\Cocont(\Q)$ is precisely the category 
of algebras for the induced KZ doctrine on $\Cat\cc(\Q)$. 
\par
We can apply all this to the quantaloid $\Idm(\Q)$, and get the following result.
\begin{proposition}\label{5}
For any small quantaloid $\Q$, $\Cocont(\Idm(\Q))$ is the category 
of algebras for the ``presheaf'' KZ doctrine 
$\P\:\Cat\cc(\Idm(\Q))\to\Cat\cc(\Idm(\Q))$.
\end{proposition}
In combination with the remarks on $\Q$-orders and $\Q$-modules in section \ref{B}, we can now justify 
the title of the paper.
\begin{theorem}\label{6}
For a small quantaloid $\Q$, the diagram
$$\Mod(\Q)\simeq\Cocont(\Idm(\Q))\xymatrix{\ar@{}[r]|{\perp}\ar@/_2mm/@<-1mm>[r]_{\cal U} & \ar@/_2mm/@<-1mm>[l]_{\cal P}}\Cat\cc(\Idm(\Q))\simeq\Ord(\Q)$$
exhibits the quantaloid $\Mod(\Q)$ as (biequivalent to) the category of algebras for the ``presheaf construction'' KZ doctrine on $\Ord(\Q)$.
\end{theorem}
As an example we shall point out how the preceding theorem is a precise generalization of the well-known fact that the internal suplattices in a localic topos $\Sh(\Omega)$ are exactly the $\Omega$-modules [Joyal and Tierney, 1984; Pitts, 1988].
\begin{example}\label{ex}
Let $\Omega$ be a locale and $(F,\leq)$ and ordered object in $\Sh(\Omega)$. We can associate to this ordered sheaf a category $\bbA$ enriched in the quantaloid $\Idm(\Omega)$ (the split-idempotent completion of the monoid $(\Omega,\wedge,1)$) as follows:
\begin{itemize}
\item objects: $\bbA_0=\coprod_{v\in\Omega} F(v)$, with types $tx=v\Leftrightarrow x\in F(v)$,
\item hom-arrows: for $x,y\in\bbA_0$, $\bbA(y,x)=\bigvee\{w\leq tx\wedge ty\mid y_{|w}\leq_w x_{|w}\}.$
\end{itemize}
That is to say, we can read off that
$$\bbA(y,x)=\mbox{``the greatest level at which $y\leq x$ in $F$''}.$$
With a slight adaptation of the arguments in [Walters, 1981; Borceux and Cruciani, 1998] one can prove that this construction extends to a (bi)equivalence of locally ordered categories $\Ord(\Sh(\Omega))\simeq\Cat\cc(\Idm(\Omega))$; the details are in [Stubbe, 2005b].
\par
We shall now explain that, under the identification of $(F,\leq)$ in $\Sh(\Omega)$ with $\bbA$, there is a {\em bijective correspondence between downsets of $F$ and presheaves on $\bbA$}; in particular do {\em principal downsets correspond with representable presheaves}.
\par
A {\em downset $S$ of $(F,\leq)$} is an $S\in\Omega^F$ (i.e.\ an $S\subseteq F_u\subseteq F$ for some $u\in\Omega$)\footnote{We write $F_u$ for the ``truncation of $F$ at $u$'' [Borceux, 1994, vol.~3, 5.2.3]: it is  the sheaf defined by $F_u(v)=F(v)$ whenever $v\leq u$ and otherwise $F_u(v)=\emptyset$.} such that
\begin{equation}\label{ex1}
(y\leq x)\wedge(x\in S)\Rightarrow(y\in S),
\end{equation}
this definition being written in the internal logic of $\Sh(\Omega)$. On the other hand, a {\em presheaf} on the $\Idm(\Omega)$-enriched category $\bbA$ is by definition a distributor $\phi\:*_{u}\dist\bbA$ for some $u\in\Idm(\Omega)$; equivalently, such is a map $\phi\:\bbA_0\to\Omega$ such that for all $x,y\in\bbA_0$, $\phi(x)\leq u\wedge tx$ and 
\begin{equation}\label{ex2}
\bbA(y,x)\wedge\phi(x)\leq\phi(y).
\end{equation}
\par
The similarity between the formulas in \eqref{ex1} and \eqref{ex2} suggests that a downset $S$ of $(F,\leq)$ is related with a presheaf $\phi$ on $\bbA$ by the clause
$$\phi(x)=\mbox{``the greatest level at which $x\in S$''.}$$
Here is how this can be made precise: Given a downset $S\subseteq F_u\subseteq F$ with its characteristic map $\varphi\:F\to\Omega$, consider the family of its components $\varphi_v\:F(v)\to\Omega(v)$ (indexed by $v\in\Omega$), extend their codomains in the obvious way to the whole of $\Omega$ and call these new maps $\phi_v\:F(v)\to\Omega$. The coproduct $\phi=\coprod_{v\in\Omega}\phi_v\:\bbA_0\to\Omega$ satisfies, for $x\in\bbA_0$,
$$\phi(x)= \bigvee\{v\leq tx\mid x_{|v}\in S(v)\}$$
so that quite obviously $\phi(x)\leq tx\wedge u$, and moreover \eqref{ex2} holds because it is just a rephrasing of \eqref{ex1}. Hence $\phi$ gives the elements of a presheaf $\phi\:*_u\dist\bbA$. Conversely, given a presheaf $\phi\:*_u\dist\bbA$ we decompose the map $\phi\:\bbA_0\to\Omega$ into a family of maps $\phi_v\:F(v)\to\Omega\:x\mapsto\phi(x)$ indexed by $v\in\Omega$. Since $\phi(x)\leq tx\wedge u$ we can restrict the codomains of each of these maps to obtain a new family
$$\left(\varphi_v\:F(v)\to\Omega(v)\:x\mapsto\phi(x)\right)_{v\in\Omega}.$$
This family is natural in $v$: Let $w\leq v$ and take any $x\in F(v)$. Then $w=\bbA(x_{|w},x)$ and therefore $w\wedge\varphi_v(x)=\bbA(x_{|w},x)\wedge\phi(x)\leq\phi(x_{|w})=\varphi_w(x_{|w})$ by \eqref{ex2}. But also $w=\bbA(x,x_{|w})$ and so, again by \eqref{ex2}, $\varphi_w(x_{|w})=\phi(x_{|w})=w\wedge\phi(x_{|w})=\bbA(x,x_{|w})\wedge\phi(x_{|w})\leq\phi(x)=\varphi_v(x)$. Thus indeed $w\wedge\varphi_v(x)=\varphi_w(x_{|w})$. Now we let $S\in\Omega^F$ be the $S\subseteq F$ with characteristic map $\varphi\:F\to\Omega$: then actually $S\subseteq F_u\subseteq F$ because $\phi(x)\leq u$, and moreover $S$ is a downset because \eqref{ex1} follows from \eqref{ex2}. The constructions $S\mapsto\phi$ and $\phi\mapsto S$ are inverse to each other under the identification of the ordered sheaf $(F,\leq)$ with the enriched category $\bbA$.
\par
In particular, the {\em principal downset $S_x$ of $F$ at $x\in F$} is the $S_x\in\Omega^F$ such that
$$(y\leq x)\Leftrightarrow (y\in S_x).$$
(Clearly such an $S_x$ is always a downset.) The corresponding presheaf $\phi_x\:*_{u}\dist\bbA$ must thus satisfy
$$\bbA(y,x)=\phi_x(y),$$
that is to say, it is the {\em representable presheaf} $\bbA(-,x)$.
\par
Now we can understand why {\em an ordered sheaf $(F,\leq)$ is an internal suplattice in $\Sh(\Omega)$ if and only if the associated $\Idm(\Omega)$-category $\bbA$ is cocomplete}: $(F,\leq)$ is an internal suplattice in $\Sh(\Omega)$ if and only if the ``principal downset inclusion'' $F\to\Omega^F$ has a left adjoint [Mikkelsen, 1976; Johnstone, 2002, B2.3.9]. But this is constructively equivalent with the existence of a left adjoint to its factorization over the (object of) downsets of $F$. By the above we know that this is the case if and only if the Yoneda embedding $Y_{\bbA}\:\bbA\mapsto\P\bbA$ has a left adjoint, which in turn means precisely that $\bbA$ is cocomplete.
\par
By \ref{6} we thus get an independent proof of the fact that the internal suplattices in $\Sh(\Omega)$ are precisely the modules on $\Omega$: $\Sup(\Sh(\Omega))\simeq\Mod(\Omega)$.  
\end{example}

\section{Module equivalence compared with sheaf equivalence}\label{D}
For any quantaloid $\Q$, let $\Z(\Q)$ be shorthand for 
$\QUANT(\Q,\Q)(\Id_{\Q},\Id_{\Q})$ and call it the {\em centre of $\Q$}. 
This is by definition a commutative quantale: that $\Z(\Q)$ 
is a quantale, is because it is an endo-hom-object of the quantaloid 
$\QUANT(\Q,\Q)$; that it is moreover commutative, is because 
$\QUANT(\Q,\Q)$ is monoidal with $\Id_{\Q}$ the unit object for the 
tensor (which is composition). Unraveling the definition, an element $\alpha\in\Z(\Q)$ is a 
family of endo-arrows
$$\Big(\endoar{A}{\alpha_A}~\Big{|}~A\in\Q_0\Big)$$
such that for every $f\:A\to B$ in $\Q$, $\alpha_B\circ f = f\circ 
\alpha_A$. Inspired by [Bass, 1968, p.~56] it is then straightforward to prove the following proposition. (Since I believe that this is a ``folk theorem'' -- and moreover the case for quantales is already mentioned in [Borceux and Vitale, 1992] -- I shall only sketch the proof.)
\begin{proposition}\label{8}
For any quantaloid $\Q$, $\Z(\Q)\cong\Z(\Mod(\Q))$. Therefore Morita-equivalent quantaloids 
have isomorphic centres.
\end{proposition}
\sketchofproof Given a natural transformation 
$\alpha\:\Id_{\Q}\to\Id_{\Q}$, build the natural transformation 
$\widehat{\alpha}\:\Id_{\Mod(\Q)}\to\Id_{\Mod(\Q)}$ whose component 
at $M\in\Mod(\Q)$ is the natural transformation 
$\widehat{\alpha}_M\:M\to M$, whose component at $A\in\Q$ is the 
$\Sup$-arrow
$$\widehat{\alpha}_M^A=M(\alpha_A)\:M(A)\to M(A).$$
Conversely, given a natural transformation 
$\beta\:\Id_{\Mod(\Q)}\to\Id_{\Mod(\Q)}$, build the natural 
transformation $\overline{\beta}\:\Id_{\Q}\to\Id_{\Q}$ whose 
component at $A\in\Q$ is the $\Q$-morphism
$$\overline{\beta}_A=\beta_{\Q(A,-)}^A(1_A)\:A\to A.$$
The mappings $\Z(\Q)\to\Z(\Mod(\Q))\:\alpha\mapsto\widehat{\alpha}$ and 
$\Z(\Mod(\Q))\to\Z(\Q)\:\beta\mapsto\overline{\beta}$ thus defined are quantale 
homomorphisms
which are each other's inverse.
\endofproof
The following is now an easy consequence.
\begin{proposition}\label{10}
For small quantaloids $\Q$ and $\Q'$,
$$\Q\simeq\Q'\impl\Ord(\Q)\simeq\Ord(\Q')\impl\Mod(\Q)\simeq\Mod(\Q')\impl\Z(\Q)\cong\Z(\Q').$$
\end{proposition}
\proof The first implication holds because ``equivalent bases 
give equivalent enriched structures''. The second implication is 
due to the monadicity explained in \ref{6}. For the third implication, see 
\ref{8}. 
\endofproof
It is an interesting problem to study the converse implications in 
the above proposition, for they do not hold in general.
However, since a quantale is commutative if and only if it equals its centre, we do have the following special case which is a refinement of the conclusion of [Borceux and Vitale, 1992].
\begin{corollary}\label{11}
For commutative quantales $\Q$ and $\Q'$,
$$\Q\simeq\Q'\iff\Ord(\Q)\simeq\Ord(\Q')\iff\Mod(\Q)\simeq\Mod(\Q').$$
\end{corollary}
A {\em locale $\Omega$} is in particular a commutative quantale, so the above 
applies. Moreover, and this in strong contrast with the case of 
quantaloids or even quantales, besides the category $\Ord(\Omega)$ of ordered sheaves 
and its subcategory $\Mod(\Omega)$ of modules (i.e.\ cocompletely ordered sheaves) on $\Omega$, we may now
also consider the category $\Sh(\Omega)$ of all sheaves. But a
locale $\Omega$ is (isomorphic to) the locale of subobjects of the 
terminal object in $\Sh(\Omega)$ (see [Borceux, 1994, vol.~3, 2.2.16] for 
example), thus we may end with the following.
\begin{corollary}\label{12}
For locales $\Omega$ and $\Omega'$,
$$\Omega\simeq \Omega'\iff\Sh(\Omega)\simeq\Sh(\Omega')\iff\Ord(\Omega)\simeq\Ord(\Omega')\iff\Mod(\Omega)\simeq\Mod(\Omega').$$
\end{corollary}

\end{document}